%
%
%
%
%
\font\tenbf=cmbx10

\font\eightrm=cmr8
\font\eightit=cmti8
\font\germ=eufm10

\def\s{\hbox{\germ S}}
\def\sectiontitle#1\par{\vskip0pt plus.1\vsize\penalty-250
\vskip0pt plus-.1\vsize\bigskip\vskip\parskip
\message{#1}\leftline{\tenbf#1}\nobreak\vglue 5pt}
\def\wt{\widetilde}
\def\wh{\widehat}
\def\ds{\displaystyle}
\def\eno{\eqalignno}
\def\ld{\lambda}

\def\al{\alpha}

\def\mod{\hbox{\rm mod~}}

\def\frac#1#2{{#1\over#2}}

\def\m@th{\mathsurround=0pt}
\def\Box{\fsquare(0.25cm, )\,}
\def\fsquare(#1,#2){
\hbox{\vrule$\hskip-0.4pt\vcenter to #1{\normalbaselines\m@th
\hrule\vfil\hbox to #1{\hfill$\scriptstyle #2$\hfill}\vfil\hrule}$\hskip-0.4pt
\vrule}}
\magnification=\magstep1
\vsize=8.6truein 
\parindent=15pt
\nopagenumbers
q-alg/9703047
\hfil
\vglue 4pc
\baselineskip=13pt
\headline{\ifnum\pageno=1\hfil\else%
{\ifodd\pageno\rightheadline \else \leftheadline\fi}\fi}
\def\rightheadline{\hfil\eightit 
Cauchy identities for universal Schubert polynomials
\quad\eightrm\folio}
\def\leftheadline{\eightrm\folio\quad 
\eightit 
Anatol N. Kirillov 
\hfil}
\voffset=2\baselineskip 
\centerline{\tenbf CAUCHY \hskip 0.1cm IDENTITIES }
\vskip 0.1cm 
\centerline{\tenbf FOR \hskip 0.1cm UNIVERSAL \hskip 0.1cm SCHUBERT \hskip 0.1cm
POLYNOMIALS
}
\vglue 16pt
\centerline{\eightrm 
ANATOL N. KIRILLOV
}
\baselineskip=18pt
\centerline{\eightit 
CRM, University of Montreal
}
\baselineskip=10pt
\centerline{\eightit 
C.P. 6128, Succursale A, Montreal (Quebec) H3C 3J7, Canada
}
\baselineskip=12pt
\centerline{\eightit
and }
\baselineskip=12pt
\centerline{\eightit 
Steklov Mathematical Institute,
}
\baselineskip=10pt
\centerline{\eightit 
Fontanka 27, St.Petersburg, 191011, Russia
}
\vglue 20pt

\centerline{\eightrm ABSTRACT}
{\rightskip=1.5pc
\leftskip=1.5pc
\eightrm\parindent=1pc
We prove the Cauchy type identities for the universal double Schubert 
polynomials, introduced recently by W.~Fulton. As a corollary, the 
determinantal formulae for some specializations of the universal double 
Schubert polynomials corresponding to the Grassmannian permutations are obtained. 
We also introduce and study the universal Schur functions and multiparameter 
deformation of Schubert polynomials. 
} 
\vglue12pt
\baselineskip=13pt
\overfullrule=1pt
\def\qed{\hfill$\vrule height 2.5mm width 2.5mm depth 0mm$}

{\bf \S 0. Introduction.}
\vskip 0.3cm

The aim of this note is to study the algebraic properties of universal 
double Schubert polynomials introduced recently by W.~Fulton [F]. More 
specifically, the main goal of our note is to generalize the results from 
[K2] on the case of universal double Schubert polynomials and universal 
Schur functions. Another goal of this note is to introduce and study the
multiparameter deformation of classical and quantum Schubert polynomials,
which (conjecturally) should be closely related with the second form of
universal double Schubert polynomials introduced by W.~Fulton [F].
Our main tool is the various generalizations of the quantum Cauchy 
identity [KM]. Based on a generalized Cauchy's identities technique, we 
give an affirmative answer on some questions raised in the paper [F]. As
another application, we prove the determinantal formulae for some
specializations of the universal double Schubert polynomials corresponding
to the Grassmannian permutations.

\vskip 0.3cm
{\bf Acknowledgments.} I would like to thank W.~Fulton for sending me his
preprint ``Universal Schubert Polynomials'' [F] which initiated this note, 
and Dr. N.A.~Liskova for fruitful 
discussions and the inestimable help with different kind of computations. 
I would like also to thank all my colleagues from CRM, University of 
Montreal and especially L.~Vinet and J.~Hanard for discussion and support.

\vskip 0.5cm
{\bf \S 1. Schur functions.}
\vskip 0.3cm

In this section we give a brief review of some basic definitions and 
results of the theory of Schur functions. In exposition we follow to the 
I.~Macdonald book [M1] (see also [M2]) where proofs and more details can 
be found.

\vskip 0.3cm
$\bullet$ Schur functions.
\vskip 0.2cm

Let us recall at first a definition of the Schur function $s_{\ld}(X_n)$ 
as the quotient of two alternates:
$$s_{\ld}(X_n)=\ds{\det\left( x_i^{\ld_j+n-j}\right)_{1\le i,j\le n}\over
\det\left( x_i^{n-j}\right)_{1\le i,j\le n}}. \eqno (1.1)
$$
where $X_n=(x_1,\ldots ,x_n)$ is the set of independent variables and 
$\ld=(\ld_1, \ldots ,\ld_n)$ is a partition of the length $\le n$.
The denominator on the RHS(1.1) is the Vandermond determinant, and is 
equal to the product $\ds\prod_{1\le i<j\le n}(x_i-x_j)$.

When $\ld =(r)$, $s_{\ld}(X_n)$ is the complete symmetric function 
$h_r(X_n)$, and when $\ld =(1^r)$, $s_{\ld}(X_n)$ is the elementary 
symmetric function $e_r(X_n)$. In terms of complete symmetric functions, 
$s_{\ld}(X_n)$ is given by Jacobi--Trudi formula
$$s_{\ld}(X_n)=\det\left(h_{\ld_i-i+j}(X_n)\right)_{1\le i,j\le n}. \eqno 
(1.2)
$$
Dually, in terms of the elementary symmetric functions, $s_{\ld}(X_n)$ is 
given by the N\"agelsbach--Kostka formula
$$s_{\ld}(X_n)=\det\left( e_{\ld'_i-i+j}(X_n)\right)_{1\le i,j\le m}, 
\eqno (1.3)
$$
where $\ld'=(\ld'_1,\ldots ,\ld'_m)$ is the conjugate ([M1], p.2) of the 
partition $\ld$.

\vskip 0.2cm
{\bf Remark 1.} Using the recurrence relation for the elementary symmetric 
functions
$$e_r(X_n)=e_r(X_{n-1})+x_ne_{r-1}(X_{n-1})
$$
it is a matter of simple row transformations to show that (cf. [K2], 
(3.2) and (4.13))
$$s_{\ld}(X_n)=\det\left(h_{\ld_i-i+j}(X_{n+1-j})\right)_{1\le i,j\le n},
$$
and dually,
$$s_{\ld}(X_n)=\det\left( e_{\ld'_i-i+j}(X_{n+j-1})\right)_{1\le i,j\le m}.
$$

To complete this Section we should mention the Cauchy identity
$$\sum_{\ld}s_{\ld}(X_n)s_{\wh\ld'}(Y_m)=\prod_{1\le i\le n}
\prod_{1\le j\le m}(x_i+y_j), \eqno (1.4)
$$
summed over all partitions $\ld =(\ld_1,\ldots ,\ld_n)$ such that $\ld_1\le 
m$, where $\wh\ld =(\wh\ld_1,\ldots ,\wh\ld_n)$ is the complementary 
partition defined by $\wh\ld_i=m-\ld_{n-i+1}$, and $\wh\ld'$ is the 
conjugate of $\wh\ld$.

\vskip 0.3cm
$\bullet$ Generalized Schur functions.
\vskip 0.2cm

The generalized Schur functions were introduced by D.~Littlewood [L]. 
Namely, for any formal series $f=\ds\sum_{k\ge 0}z^ks_k$, the generalized 
skew Schur functions are defined as the minors of the Hankel matrix
$${\bf S}(f)=(s_{j-i})_{i,j\ge 0},
$$
putting $s_i=0$, if $i<0$.

More precisely, for given partitions $\ld$ and $\mu$ such that 
$\mu\subset\ld$, $l(\ld )\le n$, the generalized skew Schur function 
$s_{\ld /\mu}(f)$ is defined to be
$$s_{\ld /\mu}(f)=\det\left( s_{\ld_i-\mu_j-i+j}\right)_{1\le i,j\le n}.
$$

{\bf Example 1.} Let us take 
$$f={\ds\prod_{y\in Y}(1-zy)\over\ds\prod_{x\in X}(1-zx)}=\sum_{k\ge 
0}s_k(X-Y)z^k.
$$
Then the generalized Schur function $s_{\ld /\mu}(f):=s_{\ld /\mu}(X-Y)$ 
coincides with the so--called super--Schur function [M1], Chapter~I, \S 
3, Example~23 .

\vskip 0.5cm
{\bf \S 2. Quantum elementary and quantum complete homogeneous polynomials.}
\vskip 0.3cm

Let $X_n=(x_1,\ldots ,x_n)$ and $Y_n=(y_1,\ldots ,y_n)$ be two sets of 
variables, and\break $q=(q_1,\ldots ,q_{n-1})$ and $q'=(q'_1,\ldots ,q'_{n-1})$ 
be two sets of independent parameters. Follow to [GK], let us define the 
quantum elementary polynomials $e_i^q(X_k)$ from the decomposition of the 
Givental-Kim determinant \ \ \ $\ds\sum_{i=0}^ke_i^q(X_k)t^{k-i}=$
\vskip 0.1cm
$$ =\det\pmatrix{ x_1+t & q_1&0 &\ldots &\ldots 
&\ldots &0\cr
-1 & x_2+t & q_2 &0 &\ldots &\ldots & 0\cr
0 & -1 &x_3+t & q_3 & 0 &\ldots & 0 \cr
\vdots &\ddots &\ddots & \ddots &\ddots &\ddots &\vdots \cr
0&\ldots & 0 &-1&x_{k-2}+t &q_{k-2} & 0 \cr
0 &\ldots &\ldots & 0 &-1 & x_{k-1}+t & q_{k-1}\cr
0 &  \ldots &\ldots &\ldots & 0 & -1 & x_k+t} \eqno (2.1)
$$
\vskip 0.4cm

We define (see [FGP], [KM], or [K2], (3.4)) the quantum complete 
homogeneous polynomial $h_k^q(X_r)$ by the following formula
$$h_k^q(X_r)=\det\left( e_{1-i+j}^q(X_{r-1+j})\right)_{1\le i,j\le k}. 
\eqno (2.2)
$$
From the very definition (2.2) of quantum complete homogeneous 
polynomials, one can deduce, [K2], (3.4), the following "inversion formula":
$$e_k^q(X_r)=\det\left( h^q_{1-i+j}(X_{r+1-j})\right)_{1\le i,j\le k}. 
\eqno (2.3)
$$
Using the recurrence relations ($0\le k\le m\le n$) for quantum elementary 
polynomials $e_k^q(X_m)$:
$$e_k^q(X_m)=e_k^q(X_{m-1})+x_me_{k-1}^q(X_{m-1})+q_{m-1}e_{k-2}^q(X_{m-2})
$$
it is a matter of simple row transformations to show that
$$h_k^q(X_r)=\det \left( e^q_{1-i+j}(X_{\min (r-1+j,n-1)})\right). \eqno 
(2.4)
$$

Now let us define the quantum super elementary polynomials 
$e_m^{q,q'}(X_k-Y_l)$, and quantum super complete homogeneous polynomials 
$h_m^{q,q'}(X_k-Y_l)$. Namely, let us put
$$\eno{
h_m^{q,q'}(X_k-Y_l)&=\sum_{j=0}^mh^q_{m-j}(X_k)e_j^{q'}(Y_l), & (2.5)\cr
e_m^{q,q'}(X_k-Y_l)&=\sum_{j=0}^me^q_{m-j}(X_k)h_j^{q'}(Y_l), & (2.6)}
$$
It follows from (2.5) and (2.6) that the quantum super elementary and 
quantum super complete 
homogeneous polynomials satisfy  the following duality
$$h_m^{q,q'}(X_k-Y_l)=e_m^{q',q}(Y_l-X_k). \eqno (2.7)
$$

More generally, let us introduce the "semi--universal" (cf. [F]) 
analogues of polynomials $h_m^{q,q'}(X_k-Y_l)$ and $e_m^{q,q'}(X_k-Y_l)$.
To do this, let us consider two sets of independent variables $c_i(j)$ 
and $d_i(j)$ for $1\le i\le j\le n-1$. It is convenient to put  
$c_i(j)=d_i(j)=1$ if $i=0$, and $d_i(j)=c_i(j)=0$ if $i<0$ or $i>j$. 
Then we define
$$\eno{
&e_m^{q'}(k|Y_l)=\sum_{j=0}^mc_{m-j}(k)h_j^{q'}(Y_l), & (2.8)\cr
&h_m^q(X_r|l)=\sum_{j=0}^mh_{m-j}^q(X_r)d_j(l). & (2.9)}
$$

Let us remark that polynomials $e_m^{q'}(k|Y_l)$ (resp. $h_m^q(X_r|k)$) 
are specialized to the quantum super elementary polynomials 
$e_m^{q,q'}(X_k-Y_l)$ (resp. to the quantum super complete homogeneous 
polynomials $h_m^{q,q'}(X_r-Y_l)$) when each $c_i(j)$ is sent to 
$e_i^q(X_j)$ (resp. each $d_i(j)$ is sent to $e_i^{q'}(Y_j)$). On the 
other hand, if in (2.8) we take $q'=0$, then 
$e_m^{q'}(k|-Y_l)|_{q'=0}=f_m(k,l,0)$, where polynomials $f_m(k,a,b)$ 
were introduced in [F], (19).

Finally, let us define the quantum super multi--Schur functions $s_{\ld 
/\mu}^{q,q'}({\cal X},{\cal Y})$, where\break ${\cal X}=(X_{k_1},\ldots ,
X_{k_m})$ and ${\cal Y}=(Y_{l_1},\ldots ,Y_{l_m})$ be 
two families of flagged sets of variables, and $\ld ,\mu$ be partitions 
of length $\le m$.
\vskip 0.2cm

{\bf Definition 1.} {\it The quantum super multi--Schur function $s_{\ld 
/\mu}^{q,q'}({\cal X},{\cal Y})$ is defined to be}
$$s_{\ld /\mu}^{q,q'}({\cal X},{\cal Y})=\det\left( 
h_{\ld_i-\mu_j-i+j}^{q,q'}(X_{k_i}-Y_{l_i})\right)_{1\le i,j\le m}. \eqno 
(2.10)
$$

\vskip 0.5cm 
{\bf \S 3. Universal double Schubert polynomials.}
\vskip 0.3cm

Follow to [F], let us define at first the universal Schubert polynomials 
${\s}_w(c)$. Let $c_i(j)$, $1\le i\le k\le n$, be set of independent 
variables. It is convenient to define $c_0(j)=1$, and $c_i(j)=0$ if $i<0$ 
or $i>j$. We start with definition of the universal Schubert polynomial 
${\s}_{w_0}(c)$ corresponding to the element of maximal length $w_0\in 
S_n$:
$${\s}_{w_0}(c,y)=\prod_{i=1}^{n-1}\left(\sum_{j=0}^iy^j_{n-i}c_{i-j}(i)\right). 
\eqno (3.1)
$$
The next step is to define the double polynomials ${\s}_w(c,y)$ for any 
permutation $w\in S_n$. Let us put [F], cf [KM],
$${\s}_w(c,y)=\partial_{ww_0}^{(y)}{\s}_{w_0}(c,y), \eqno (3.2)
$$
where divided difference operator $\partial_{ww_0}^{(y)}$ acts on the $y$ 
variables.

\vskip 0.2cm
{\bf Definition 2} ([F]). {\it Let $w\in S_n$ be a permutation. The single 
universal Schubert polynomial ${\s}_w(c)$ is defined to be the 
specialization $y_1=\cdots =y_n=0$ of the double polynomial ${\s}_w(c,y)$:}
$${\s}_w(c)={\s}_w(c,0).
$$

Finally, let us define (see [F]) the universal double Schubert 
polynomials ${\s}_w(c,d)$, where $c$ stands for the variables $c_i(j)$ 
and $d$ stands for another set of variables $d_i(j)$. 

\vskip 0.2cm
{\bf Definition 3} ([F]). {\it The universal double Schubert polynomial
${\s}_w(c,d)$ is defined by the following formula
$${\s}_w(c,d)=\sum_{u\in S_n}{\s}_u(c){\s}_{uw^{-1}}(d), \eqno (3.3)
$$
where the sum is over all $u\in S_n$ such that $l(u)+l(uw^{-1})=l(w)$.}

It is clear from (3.3) that 
$${\s}_w(c,d)={\s}_{w^{-1}}(d,c). \eqno (3.4)
$$

\vskip 0.5cm
{\bf \S 4. Generalized quantum Cauchy identity.}
\vskip 0.3cm

The quantum Cauchy identity is the quantum analog of the Cauchy formula 
in the theory of Schubert polynomials ([M1], (5.10)). As it was shown in 
[KM], Section~4, the quantum Cauchy identity corresponds to the 
quantization of Cauchy's formula (1.4) with respect to the $x$ variables.

The generalized quantum Cauchy identity corresponds to the quantization 
of quantum Cauchy identity ([KM], [K1]) with respect to the $y$ variables.

\vskip 0.2cm
{\bf Theorem 1.} (Generalized quantum Cauchy identity for quantum 
Schubert polynomials)
$$\sum_{w\in S_n}{\s}_w^q(X_n){\s}_{ww_0}^{q'}(Y_n)=
{\s}_{w_0}^{q,q'}(X_n,Y_n), \eqno (4.1)
$$
{\it where ${\s}_{w_0}^{q,q'}(X_n,Y_n)$ is given by }
$${\s}_{w_0}^{q,q'}(X_n,Y_n)=\det\left( 
h_{n-2i+j}^{q,q'}(X_i-Y_{n-i})\right)_{1\le i,j\le n-1}.
$$

Theorem~1 follows from the following more general result:

\vskip 0.2cm
{\bf Theorem 2.}
$$\sum_{w\in S_n}{\s}_w^q(X_n){\s}_{ww_0}(d)={\s}_{w_0}^q(X_n,d), 
\eqno (4.2)
$$
{\it where
$${\s}_{w_0}^q(X_n,d)=\det\left( h_{n-2i+j}^q(X_i| n-i)\right)_{1\le 
i,j\le n-1},
$$
and $h_m^q(X_r|l)$ is defined by (2.9).}

{\it Proof of Theorem~2.} First of all, if all $q_i=0$, the formula (4.2) 
follows from [F], Lemma~2.1. General case follows from the following 
statement:
\vskip 0.2cm

{\bf Lemma 1.} {\it Let $c_i(j)$ and $d_i(j)$ for $1\le i\le j\le n-1$ be 
two sets of independent variables, and let
$$e_m(k|l):=\sum_{j=0}^mc_{m-j}(k)d_j(l)
$$
be "universal" elementary polynomial. Then
$$\det\left( e_{n-2i+j}(i|n-i)\right)_{1\le i,j\le 
n-1}=\sum_{I\subset\delta_n}s_Id_{\delta_n-I}, \eqno (4.3)
$$
where}
$$\eno{
&s_I:=\det\left( c_{i_{\al}-\al +\beta}(\al )\right)_{1\le\al ,\beta\le 
n-1}, \  \  if \ \ I=(i_1,\ldots ,i_{n-1});\cr
&d_J:=\prod_{k=1}^{n-1}d_{j_k}(n-k);\cr
&c_0(j)=d_0(j)=1, \ \ c_i(j)=d_i(j)=0, \ if \ i<0 \ \ or \ \ i>j.}
$$

Indeed, using Lemma~1, one can rewrite the formula (4.2) with $q=0$ in 
the following form
$$\sum_{w\in S_n}{\s}_w(X_n){\s}_{ww_0}(d)={\s}_{w_0}(X_n,d)=
\sum_{I\subset\delta_n}s_I(X_1,\ldots ,X_{n-1})d_{\delta_n-I}=
\sum_{I\subset\delta_n}x^Id_{\delta_n-I}. \eqno (4.4)
$$
We used here the following formula ([M3], (3.5'))
$$x^I=\det\left( h_{i_{\al}-\al +\beta}(X_{\al})\right)_{1\le\al ,\beta\le 
n-1}=s_I(X_1,\ldots ,X_n).
$$
Hence, using the properties of quantization map ([FGP], [KM], [K2]) with 
respect to the $x$ variables, we obtain
$$\eno{
\sum_{w\in S_n}{\s}_w^q(X_n){\s}_{ww_0}(d)&=\sum_{I\subset\delta_n}\wt 
x^Id_{\delta_n-I}=\sum_{I\subset\delta_n}s_I^q(X_1,\ldots 
,X_{n-1})d_{\delta_n-I}\cr \cr
&=\det\left( h_{n-2i+j}^q(X_i|n-i)\right)_{1\le 
i,j\le n-1}={\s}_{w_0}^q(X_n,d).}
$$
Here we used the following result ([K2], Corollary~4): if 
$I\subset\delta_n$, then $\wt x^I=s_I^q(X_1,\ldots ,X_{n-1})$, where $\wt 
x^I$ is the quantization of monomial $x^I=x_1^{i_1}\ldots x_{n-1}^{i_{n-1}}$, 
and
$$s_I^q(X_1,\ldots ,X_{n-1})=\det\left( h_{i_{\al}-\al 
+\beta}^q(X_{\al})\right)_{1\le\al ,\beta\le n-1}.
$$
\qed

\vskip 0.2cm
{\bf Corollary 1.}
$$\eno{
&~~~~~~~~~~~\sum_{w\in S_n}{\s}_w^{q,q"}(X_n,Z_n){\s}_{ww_0}^{q',q"}(Y_n,-Z_n)=
{\s}_{w_0}^{q,q'}(X_n,Y_n), & (4.5)\cr
&\sum_{\matrix{~_{u\in S_n}\cr ~_{l(u)+l(uw^{-1})=l(w)}}}
{\s}_u^{q,q"}(X_n,Z_n){\s}_{uw^{-1}}^{q',q"}(Y_n,-Z_n)=
{\s}_w^{q,q'}(X_n,Y_n). & (4.6)}
$$
Here the double quantum Schubert polynomials ${\s}_w^{q,q'}(X_n,Y_n)$ is 
defined as follows (cf. [F])
$${\s}_w^{q.q'}(X_n,Y_n)=\sum_{u\in 
S_n}{\s}_u^q(X_n){\s}_{uw^{-1}}^{q'}(Y_n),
$$
where the sum is over all $u\in S_n$ such that $l(u)+l(uw^{-1})=l(w)$.
\vskip 0.2cm

{\bf Remark 2.} It follows from (4.6) that the double quantum Schubert 
polynomial ${\s}_w^{q,q'}(X_n,Y_n)$ can be obtained from the universal 
double Schubert polynomial ${\s}_w(c,d)$ (see [F]) under the specialization
$$\eno{
&c_i(k)\longrightarrow e_i^q(X_k),\cr
&d_j(l)\longrightarrow e_j^{q'}(Y_l).}
$$

{\bf Remark 3.} Polynomial ${\s}_{w_0}^{q,q'}(X_n,Y_n)$ appears at first in 
[KM], Remark~11, and corresponds to the dual class of the diagonal in the 
quantum cohomology ring 
$$QH^*(Fl_n,q)\otimes QH^*(Fl_n,q'). 
$$
The specialization $\Phi (x):={\s}_w^{q,q}(X_n,X_n)$ plays an important role 
in the proof of the Vafa--Intriligator type formula for the quantum 
cohomology ring of flag  variety ([KM], Section~8.1). On the other hand, 
it follows from (4.2) that ${\s}_w^{q,q}(X_n,-X_n)=0$, if $w\ne {\rm id}$.

\vskip 0.5cm
{\bf \S 5. Cauchy identity for universal Schubert polynomials.}
\vskip 0.3cm

In this Section we are going to formulate the Cauchy identity for the 
universal Schubert polynomials (see [F], or Section~3). We start with 
definition of polynomials ${\s}_I(c)$,\break $I\subset\delta_n:=
(n-1,n-2,\ldots ,1,0)$. Let us remark, that if $I\subset\delta_n$, then the 
monomial $x^I:=x_1^{i_1}\cdots x_{n-1}^{i_{n-1}}$ is a ${\bf Z}$--linear 
combination of elementary polynomials
$$\eno{
&e_J(X_{n-1}):=\prod_{k=1}^{n-1}e_{j_k}(X_{n-k}), \ \ {\rm say,}\cr
&x^I=\sum_{J\subset\delta_n}\al_{I,J}e_J(X_{n-1}),& (5.1)}
$$
and such representation (5.1) is unique ([LS]).

Let us define 
$${\s}_I(c)=\ds\sum_{J\subset\delta_n}\al_{I,J}c_J, \eqno (5.2)
$$
where $c_J=\ds\prod_{k=1}^{n-1}c_{j_k}(n-k)$.

\vskip 0.2cm
{\bf Theorem 3} (Cauchy identity for universal Schubert polynomials).
$${\s}_{w_0}(c,d)=\sum_{I\subset\delta_n}{\s}_I(c)d_{\delta_n-I}, \eqno 
(5.3)
$$
{\it where $d_J=\ds\prod_{k=1}^{n-1}d_{j_k}(n-k)$.}

{\it Proof.} Under the specialization $c_i(k)\to e_i^q(X_k)$, the formula 
(5.3) is reduced to (4.2). General case follows from [F], Lemma~2.1.
\qed

\vskip 0.3cm

{\bf Corollary 2.} {\it Let $b=\{ b_j(k)\}$ be the third set of variables 
such that $b_0(k)=1$, and $b_j(k)=0$, if $j<0$ or $j>k$, then}
$$\eno{
&\sum_{w\in S_n}{\s}_w(c,b){\s}_{ww_0}(d,\wt b)={\s}_{w_0}(c,d), & (5.4)\cr
&\sum_{u\in S_n}{\s}_u(c,b){\s}_{uw^{-1}}(d,\wt b)={\s}_w(c,d), &(5.5)}
$$
{\it where $\wt b_j(k)=(-1)^jb_j(k)$.}

\vskip 0.2cm
{\bf Corollary 3.} {\it ${\s}_w(b,\wt b)=0$, if $w\ne{\rm id}$.}

{\it Proof.} Let us take $b=d$ in (5.5).
\qed

\vskip 0.4cm
{\bf Remark 4.} Corollary~3 was formulated by W.~Fulton as a Conjecture 
(see [F], (30)).

\vskip 0.5cm
{\bf \S 6. Cauchy identity for universal Schubert polynomials of the 
second form.}
\vskip 0.3cm

Let us remind [F] the definition of the universal Schubert polynomials of 
the second form. Follow to [F], let us consider a set of variables 
$g_i[j]$ for $i\ge 1$ and $j\ge 0$ with $i+j\le n$; and let us put ${\rm 
deg}(g_i[j])=j+1$. Now we are going to define the generalized elementary 
functions $\Box_i(k)$. Let us consider the generalized Givental--Kim 
determinant (cf. [F])
$\Box_k(t,g)=\det (tI_k+A_k)$, where $A_k$ is the $k\times k$ matrix 
with $g_i[j-i]$ in the $(i,j)$ position for $i\le j$, and with -1 in 
position $(i+1,i)$ below the diagonal. 

\vskip 0.2cm
{\bf Definition 4.} {\it The generalized 
elementary functions $\Box_i(k)$ are defined from the decomposition}
$$\Box_k(t,g)=\sum_{j=0}^kt^{k-j}\Box_j(k). \eqno (6.1)
$$

\vskip 0.2cm
{\bf Definition 5} ([F]). {\it The universal Schubert polynomial of the 
second form, denoted ${\s}_w(g)$, is obtained from ${\s}_w(c)$ by 
replacing each $c_i(k)$ by the generalized elementary function 
$\Box_i(k)$.}

Cauchy's type identity for polynomials $\s_w(g)$, can be obtained as a 
simple corollary of the Cauchy identity for 
universal Schubert polynomials (Theorem~3). 

\vskip 0.2cm
{\bf Corollary 4.} {\it Let us define \ \ 
${\s}_{w_0}(g,Y_n):=\ds\prod_{j=1}^{n-1}\Box_j(y_{n-j},g).$ \ \ 
Then
$$\eno{
&\sum_{w\in S_n}{\s}_w(g){\s}_{ww_0}(Y_n)={\s}_{w_0}(g,Y_n),& (6.2)\cr
&\sum_{w\in S_n}{\s}_w(g,Z_n){\s}_{ww_0}(Y_n,-Z_n)={\s}_{w_0}(g,Y_n),& 
(6.3)}
$$
where ${\s}_w(g,Z_n):=\partial_{ww_0}^{(z)}{\s}_{w_0}(g,Z_n)$.}

\vskip 0.2cm
{\bf Remark 5.} The classical Schubert polynomials can be recovered from 
${\s}_w(g)$ by setting $g_i[0]=x_i$ and $g_i[j]=0$ for $j\ge 1$, and the 
quantum Schubert polynomial ${\s}_w^q$ can be recovered from ${\s}_w(g)$ 
by setting $g_i[0]=x_i$, $g_i[1]=q_i$, and $g_i[j]=0$ for $j\ge 2$.

\vskip 0.5cm
{\bf \S 7. Orthogonality.}
\vskip 0.3cm

Let us put $g_i[0]=x_i$, $1\le i\le n$, and consider the polynomial ring 
${\bf Z}[g]$ in variables $g_i[j]$ for $j>0$ and $i+j\le n$. Follow [F], 
let us consider a universal ring
$$R_n={\bf Z}[g][x_1,\ldots ,x_n]/{\cal R},
$$
where the ideal ${\cal R}\subset{\bf Z}[g][x_1,\ldots ,x_n]$ is generated 
by the generalized elementary functions $\Box_i(n)$, $1\le i\le n$. Here 
we consider $\Box_i(n)$ as polynomial in $x$ and $g$ variables. There 
exists a natural pairing $\langle ,\rangle_{\cal R}$ on the ring $R_n$ 
with values in ${\bf Z}[g]$, which is induced by the Grothendieck residue 
with respect to the ideal ${\cal R}$. Namely, if $f\in{\bf Z}[g][X_n]$, 
then consider the image $\overline f$ of $f$ in the quotient ring $R_n$, 
and then picking off the coefficient of ${\s}_{w_0}(g)$, or the 
coefficient of $x_1^{n-1}x_2^{n-2}\cdots x_{n-1}$.

\vskip 0.2cm
{\bf Conjecture 1} (cf. [F], (29)). {\it Universal Schubert polynomials 
${\s}_w(g)$ can be obtained as the Gram-Schmidt orthogonalization of the 
set of lexicographically ordered monomials $\{ x^I|I\subset\delta_n\}$ 
with respect to the residue pairing $\langle ,\rangle_{\cal R}$:
$$\eno{
&\langle{\s}_u(g),{\s}_v(g)\rangle_{\cal R}=\langle{\s}_u,{\s}_v\rangle =
\cases{1, &if $v=w_0u$,\cr 0, & otherwise ;}\cr
&{\s}_w(g)=x^{c(w)}+\sum_{I<c(w)}a_I(g)x^I,}
$$
where $a_I(g)\in{\bf Z}[g]$, and $I<c(w)$ means the lexicographical order,
and $c(w)$ is the code of a permutation $w\in S_n$, [M3], p.9.}

It is clear from (6.1) that
$${\s}_{w_0}(g,Y_n)=\sum_{I\subset\delta_n}\Box_I(g)y^{\delta_n-I}, \eqno 
(7.1)
$$
where
$$\Box_I(g)=\prod_{k=1}^{n-1}\Box_{i_k}(n-k). \eqno (7.2)
$$

{\bf Conjecture 2.} {\it If $I\subset\delta_n$, and $J\subset\delta_n$, 
then
$$\langle\,\Box_I(g),\Box_J(g)\rangle_{\cal R}=\langle e_I(x),e_J(x)\rangle 
, \eqno (7.3)
$$
where $e_I(x)=\ds\prod_{k=1}^{n-1}e_{i_k}(X_{n-k})$.}

As it was shown in [KM], Conjecture~1 follows from Conjecture~2. 
Conjecture~2 can be proven on the way suggested in [KM], using the 
recurrence relations for $\Box_I(g)$. Details will appear elsewhere.

\vskip 0.5cm
{\bf \S 8. Multiparameter deformation of Schubert polynomials.}
\vskip 0.3cm

In the joint paper with S.~Fomin [FK] we introduced the 
quadratic algebra ${\cal E}_n^q$ ([FK], Definition~2.1) which is closely 
related to the small quantum cohomology ring of the flag variety. More 
precisely, we construct a commutative subalgebra in ${\cal E}_n^q$ which 
appears to be canonically isomorphic to the small quantum cohomology ring 
of the flag variety ([FK], [P]). This construction admits a natural 
generalization (see [FK], Section~15) which gives rise to another natural 
generalization of the classical and quantum Schubert polynomials. Let us 
briefly explain this construction. Namely, let us consider the set of 
variables\break $t=\{ t_{ij}|1\le i<j\le n\}$, and the polynomial ring ${\bf 
Z}[t][X_n]:={\bf Z}[t][x_1,\ldots ,x_n]$. Let ${\cal T}$ be an ideal in 
${\bf Z}[t][X_n]$ which is generated by the following generalization of 
elementary symmetric functions $e_m(X_k)$:
$$e_m(t|X_k)=\sum_l\sum_{\matrix{~_{1\le i_1<\cdots <i_l\le n}\cr ~_{j_1>i_1,
\ldots ,j_l>i_l}}}e_{m-2l}(X_{\overline{I\cup J}})\prod_{k=1}^lt_{i_kj_k},
\eqno (8.1)
$$
where $i_1,\ldots ,i_l,j_1,\ldots ,j_l$ should be distinct elements of 
the set $\{ 1,\ldots ,n\}$, and $X_{\overline{I\cup J}}$ denotes the set 
of variables $x_a$ for which the subscript $a$ is neither one of the 
$i_k$ nor one of the $j_k$.

Let us consider the quotient ring
$$B_n:={\bf Z}[t][X_n]/{\cal T}.
$$
It is easy to see that ${\rm dim}B_n=n!$, and there exists a natural 
pairing $\langle ,\rangle_{\cal T}$ on the ring $B_n$ which is given by 
the Grothendieck residue with respect to the ideal ${\cal T}$ (see 
Section~7). Follow to the general strategy of [KM], we define a new 
multiparameter deformation of the quantum Schubert polynomials:

\vskip 0.2cm
{\bf Definition 6.} {\it Define the multiparameter Schubert polynomials, 
denoted by ${\s}_w^t$, as the Gram--Schmidt orthogonalization of the set 
of lexicographically ordered monomials\break $\{ x^I|I\subset\delta_n\}$ with 
respect to the residue pairing $\langle ,\rangle_{\cal T}$:
$$\eno{
&\langle{\s}_u^t,{\s}_v^t\rangle_{\cal T}=\langle{\s}_u,{\s}_v\rangle 
=\cases{1, & if $v=w_0$,\cr 0,& otherwise;}\cr 
& \cr
&{\s}_w^t=x^{c(w)}+\sum_{I<c(w)}a_I(t)x^I,}
$$
where $a_i(t)\in{\bf Z}[t]$, and $I<c(w)$ means the lexicographical order.}

One can easily check, that if $w\in S_n$, $x^I=x_1^{i_1},\ldots 
,x_n^{i_n}$, then
$$\langle w(x^I)\rangle_{\cal T}=(-1)^{l(w)}w\left(\langle 
x^I\rangle_{\cal T}\right) .\eqno (8.2)
$$

{\bf Example 2.} For the symmetric group $S_3$ we have
$$\eno{
&\langle x_1^3x_2^2\rangle_{\cal T}=-2t_{12}-t_{13},\cr
&\langle x_1^4x_2\rangle_{\cal T}=t_{12}+2t_{13}}
$$
(Hint: $x_1^3\equiv t_{12}(2x_1+x_2)+t_{13}(x_1-x_2)~\mod{\cal T}$).

Consequently,
$$\eno{
&{\s}_{321}^t=x_1^2x_2+t_{12}x_1-t_{13}x_2,~~~~~~~~~~~~~~~~~~~~~~~~~~~\cr
&{\s}_{231}^t=x_1x_2+t_{12},\cr
&{\s}_{321}^t=x_1^2-t_{12}-t_{13},\cr
&{\s}_{132}^t=x_1+x_2,\cr
&{\s}_{213}^t=x_1,\cr
&{\s}_{123}^t=1.}
$$

{\bf Example 3.} For the symmetric group $S_4$ we have
$$\eno{
\langle x_1^5x_2^2x_3\rangle_{\cal T}&=t_{12}+2t_{13}+3t_{14},\cr
\langle x_1^4x_2^3x_3\rangle_{\cal T}&=-3t_{12}-t_{13}-2t_{14},\cr
\langle x_1^4x_2^2x_3^2\rangle_{\cal T}&=2t_{12}-2t_{13},\cr
\langle x_1^3x_2^3x_3^2\rangle_{\cal T}&=2t_{13}+t_{14}-2t_{23}-t_{24},\cr
\langle x_1^5x_2^3\rangle_{\cal T}&=t_{14}-t_{13},\cr
\langle x_1^5x_2^3x_3^2\rangle_{\cal T}&=-3t_{12}^2+2t_{13}^2+t_{14}^2+
8t_{12}t_{13}+6t_{13}t_{14}+t_{12}t_{14}-2t_{12}t_{23}\cr
&~~~~-t_{12}t_{24}
-t_{13}t_{34}-4t_{13}t_{23}-6t_{14}t_{23}-4t_{14}t_{24}+t_{14}t_{34},\cr
&~~~~{\rm and~~so~~on}.}
$$

Using the symmetry property (8.2) and residue formulae above, one can find 
all multiparameter Schubert polynomials ${\s}_w^t$, $w\in S_4$. We will 
give the answer only for ${\s}_{w_0}^t$. Namely,
$$\eno{
{\s}_{4321}^t&=x_1^3x_2^2x_3+2t_{12}x_1^2x_2x_3-(t_{13}+2t_{14})x_1x_2^2x_3+
t_{13}x_1^2x_2^2\cr
&+(t_{14}+t_{23})x_1^3x_2-t_{24}x_1^3x_3+
(t_{12}t_{14}+t_{12}t_{23}-t_{13}t_{24})x_1^2\cr
&+(t_{12}t_{13}-t_{12}t_{14}-2t_{13}t_{14}-t_{13}t_{23}
-t_{14}^2-2t_{14}t_{23})x_1x_2\cr
&+(t_{12}^2+t_{12}t_{14}+t_{12}t_{24}+t_{13}t_{24}+t_{14}t_{24})x_1x_3\cr
&+(-t_{13}^2+t_{14}^2-t_{13}t_{14}-t_{13}t_{24}+t_{14}t_{34})x_2^2\cr
&+(-t_{12}t_{13}-t_{12}t_{14}-t_{13}t_{14}+t_{14}^2)x_2x_3.}
$$
There exists an alternative way to compute the multiparameter Schubert 
polynomials, using the algebra ${\cal E}_n^t$ from [FK], Section~15. We 
are going to present detailed exposition in a separate publication. Let 
us say only, that the multiparameter Schubert polynomials have many nice 
properties such as stability, orthogonality, Pieri's type formula (cf. 
[FK], [P]), and so on. However, the structural constants $\al_{u,v}^w(t)$:
$${\s}_u^t\cdot{\s}_v^t=\sum_w\al_{u,v}^w(t){\s}_w^t
$$
appear to be the polynomials in $t$'s which may have, in general, some
negative coefficients.

It seems an interesting task to understand a connection between universal 
Schubert polynomials ${\s}_w(g)$ [F], and multiparameter Schubert 
polynomials ${\s}_w^t$.

\vskip 0.5cm
{\bf \S 9. Universal Schur and universal factorial Schur functions.}
\vskip 0.3cm

In this section, we introduce the universal Schur functions $s_{\ld}(h)$, 
and universal factorial Schur functions $s_{\ld}(g,h)$, and study their 
basic properties.

\vskip 0.2cm
{\bf Definition 7.} {\it Let $\ld$ be a partition, $\ld\subset 
((n-r)^r)$, $1\le r<n$. The universal Schur function $s_{\ld}(g)$ is 
defined as the universal Schubert polynomial ${\s}_w(g)$, corresponding 
to the Grassmannian permutation $w\in S_n$ of shape $\ld$ and descent at 
$r$.}

\vskip 0.2cm
{\bf Definition 8.} {\it Let $\ld$ be a partition such that $\ld_1\le 
n-r$ and $l(\ld )\le r$ for some $r$. We define a universal factorial 
Schur function $s_{\ld}(g,h)$ to be equal to the universal double 
Schubert polynomial ${\s}_w(g,h)$, where $w\in S_n$ is the Grassmannian 
permutation of shape $\ld$ and descent at $r$.}

Now we are going to explain a connection between the universal Schur 
functions and the Macdonald construction of the "9-th variation" of Schur 
functions, [M1], Chapter~1, \S 3, Example~21, or [M2], 9-th Variation. Let 
us remind at first some definitions from [M2].

Let $h_i(k)$ ($i\ge 1$, $k\in{\bf Z}$) be independent 
indeterminates over ${\bf Z}$. Also, for convenience, define $h_0(k)=1$ 
and $h_i(k)=0$ for $i<0$ and all $k\in {\bf Z}$. Define an automorphism 
of the ring $R$ generated by the $h_i(k)$ by
$$\varphi (h_i(k))=h_i(k+1), \ \ {\rm for~all} \ \ i,k.
$$
Thus $h_i(k)=\varphi^kh_i$, where $h_i:=h_i(0)$.

Now define [M2], for any two partitions $\ld ,\mu$ of length $\le n$,
$$s_{\ld /\mu}(h)=\det\left( 
\varphi^{\mu_j-j+1}h_{\ld_i-\mu_j-i+j}\right)_{1\le i,j\le n}, \eqno (9.1)
$$
and in particular ($\mu =0$)
$$s_{\ld}(h)=\det\left( \varphi^{-j+1}h_{\ld_i-i+j}\right)_{1\le i,j\le 
n}. \eqno (9.2)
$$

From (9.1) it follows that $h_r=s_{(r)}$ ($r\ge 0$), and we define
$$e_r:=e_r(h)=s_{(1^r)}(h), \eqno (9.3)
$$
for all $r\ge 0$, and $e_r=0$ for $r<0$.

Schur functions $s_{\ld}:=s_{\ld}(h)$ possess many properties similar to 
those of the classical Schur functions. For example,

$\bullet$ $s_{\ld /\mu}=0$, unless $\ld\supset\mu$;

\vskip 0.1cm
$\bullet$ N\"agelsbach--Kostka's formula:
$$s_{\ld /\mu}=\det\left(\varphi^{-\mu_j'+j-1}e_{\ld_i'-\mu_j'-i+j}
\right)_{1\le i,j\le n}; \eqno (9.4)
$$

$\bullet$ Giambelli's formula: if $\ld =(\al_1,\ldots 
,\al_r|\beta_1,\ldots ,\beta_r)$ in Frobenius notation, then
$$s_{\ld}=\det\left( s_{(\al_i|\beta_j)}\right)_{1\le i,j\le r}. \eqno 
(9.5)
$$

See [M1], [M2] where the proofs and additional results are contained. For 
quantum Schur functions the formulae (9.4) and (9.5) were proven in [K2].

Now we are ready to explain a connection between the universal Schur 
functions and the Macdonald 9-th variation of Schur functions.

\vskip 0.2cm
{\bf Proposition 1.} {\it Let $w\in S_n$ be a Grassmannian permutation 
with shape $\ld$ and descent at $r$. Then
$${\s}_w(c)=\varphi^rs_{\ld}(c), \eqno (9.6)
$$
where the Schur function $s_{\ld}(c)$ is obtained from $s_{\ld}(h)$ (see  
(9.1)) by setting $h_i(j)\to c_i(j)$.}

{\it Proof} follows from [F], Proposition~4.4, or Section~10.

It follows from Proposition~1, that the universal Schur functions 
$s_{\ld}(g)$ satisfy (9.4) and (9.5) with $\mu =0$ (cf. [K2]).

\vskip 0.5cm
{\bf \S 10. Determinantal formulae.}
\vskip 0.3cm

In this Section we give a generalization of some determinantal formulae 
obtained in [K2], Section~4, and [F], Proposition~4.4. To start, let us 
remind a few definitions:

$\bullet$ ${\s}_w^q(X_n,d)$ is the specialization $c_j(k)\to e_j^q(X_k)$ 
of the universal double Schubert polynomial ${\s}_w(c,d)$;

$\bullet$ ${\s}_w^{q'}(c,Y_n)$ is the specialization $d_j(k)\to 
e_j^{q'}(Y_k)$ of ${\s}_w(c,d)$;

$\bullet$ if $q'=0$ then ${\s}_w^{q'}(c,Y_n)|_{q'=0}={\s}_w(c,y)$, where 
double polynomial ${\s}_w(c,y)$ is defined in [F], or Section~3, (3.2).

\vskip 0.2cm
{\bf Theorem 4.} {\it Let $w\in S_n$ be a Grassmannian permutation with 
shape $\ld$ and descent at $r$, and let $\wt\theta =\theta (w^{-1})$ be 
the flag of inverse permutation $w^{-1}$. Then}
$${\s}_w^{q'}(c,Y_n)=\det\left( e^{q'}_{\ld_i'-i+j}
(r-1+j|Y_{\wt\theta_i})\right)_{1\le i,j\le n-r}. \eqno (10.1)
$$

{\it Proof.} For $q'=0$, the formula (10.1) coincides with 
Proposition~4.4 in [F]. For general $q'$, (10.1) can be obtained from the
particular case $q'=0$ by applying the quantization map with respect to 
the $y$ variables.

\qed

\vfil\eject
{\bf References.}
\vskip 0.5cm

\item{[F]} Fulton W., {\it Universal Schubert polynomials,} Preprint, 1997,
      alg-geom/9702012, 20p.;

\item{[FGP]}  Fomin S., Gelfand S. and Postnikov A., {\it Quantum Schubert
     polynomials,} Preprint, 1996, \hbox{AMSPPS} \#199605--14--008, 44p.;

\item{[FK]} Fomin S. and Kirillov A.N., {\it Quadratic algebras, Dunkl 
      elements, and Schubert calculus,} Preprint, 1997, 
      AMSPPS \#199703--05--001, 34p.;
        
\item{[GK]} Givental A. and Kim B., {\it Quantum cohomology of flag manifolds and 
        Toda lattices}, Comm. Math. Phys., 1995, v.168, p.609-641;

\item{[K1]} Kirillov A.N., {\it Quantum Grothendieck polynomials,} 
          Preprint, 1996, q-alg/9610034, 12p;        

\item{[K2]} Kirillov A.N., {\it Quantum Schubert polynomials and quantum 
      Schur functions,} Preprint, 1997, CRM-2452, 20p. {\it and} q-alg/9701005;

\item{[KM]} Kirillov A.N. and Maeno T., {\it Quantum double Schubert 
 polynomials, quantum Schubert polynomials and Vafa--Intriligator formula,}
 Preprint, 1996, q-alg/9610022, 52p.;

\item{[L]} Littlewood D.E., {\it The theory of group characters}, 2nd 
ed., Oxford University Press, 1950;

\item{[LS]} Lascoux A. and Sch\"utzenberger M.-P., {\it Functorialit\'e 
des polynomes de Schubert,} Contemporary Math., 1989, v.88, p.585-598;

\item{[M1]} Macdonald I., {\it Symmetric functions and Hall 
polynomials,} Second ed., Oxford Univ. Press, New York/London, 1995;

\item{[M2]} Macdonald I., {\it Schur functions: theme and variations,} 
Publ. I.R.M.A. Strasbourg, 1992, 498/S-28, Actes 28-e Seminaire Lotharingien, 
p.5--39;

\item{[M3]} Macdonald I., {\it Notes on Schubert polynomials}, Publ.
 LCIM, 1991, Univ. of Quebec a Montreal;

\item{[P]} Postnikov A., {\it On quantum Monk's and Pieri's formulas,} 
      Preprint, 1997, AMSPPS \#199703--05--002, 12p.

\end